\documentclass[12pt]{article}
\usepackage{amsmath}
\usepackage{amssymb}
\usepackage{amsthm}
\usepackage{graphicx}
\usepackage{stmaryrd}
\usepackage{pst-all,pst-node,pstricks,pst-plot,pst-3dplot}
\newtheorem{thm}{Theorem}
\newtheorem{lem}{Lemma}
\newtheorem{cor}{Corollary}

\begin{document}

\begin{center}
{\bf\LARGE  A de Bruijn - Erd\H{o}s theorem and metric spaces}\\
\vspace{1cm}
Ehsan Chiniforooshan and Va\v{s}ek Chv\'{a}tal 
\\
{\em
Department of Computer Science and Software Engineering\\
Concordia University\\
Montreal, Quebec H3G 1M8, Canada\\
}
\end{center}

\vspace{0.7cm}

\begin{center}
  {\bf Abstract}
\end{center}
  
\begin{center}
\begin{minipage}{12cm}
 {\small De Bruijn and Erd\H os proved that every
    noncollinear set of $n$ points in the plane determines at least
    $n$ distinct lines. Chen and Chv\'{a}tal suggested a possible
    generalization of this theorem in the framework of metric
    spaces. We provide partial results in this direction.  }
\end{minipage}
\end{center}
\vspace{0.3cm}

\section{Introduction}

Two distinct theorems are referred to as ``the de Bruijn - Erd\H os
theorem''. One of them \cite{DE51} concerns the chromatic number of
infinite graphs; the other \cite{DE48} is our starting point: 
\begin{center}
{\em
Every noncollinear set of $n$ points in the plane\\ 
determines at least $n$ distinct lines.\/}
\end{center}
This theorem involves neither measurement of distances nor measurement
of angles: the only notion employed here is incidence of points and
lines. Such theorems are a part of {\em ordered geometry\/}
\cite{C61}, which is built around the ternary relation of {\em
betweenness\/}: point $y$ is said to lie between points $x$ and $z$ if
$y$ is an interior point of the line segment with endpoints $x$ and
$z$. It is customary to write $[xyz]$ for the statement that $y$ lies
between $x$ and $z$. In this notation, a {\em line\/} $\overline{uv}$
is defined \ ---\  for any two distinct points $u$ and $v$ \ ---\  as
\begin{equation}\label{def.first}
\{p: [puv]\}\;\cup\; \{u\}\;\cup\; \{p: [upv]\}\;\cup\; \{v\}\;\cup\; \{p: [uvp]\}. 
\end{equation}

In terms of the Euclidean metric $d$, we have
\begin{equation}\label{def.btw}
\mbox{ $[abc] \;\Leftrightarrow\; a,b,c$ are three distinct points and
$d(a,b)+d(b,c) = d(a,c)$.}
\end{equation}
For an arbitrary metric space, equivalence (\ref{def.btw}) defines the
ternary relation of {\em metric betweenness} introduced in \cite{Men}
and further studied in \cite{Blu,Bus,Chv}; in turn, (\ref{def.first})
defines the line $\overline{uv}$ for any two distinct points $u$ and
$v$ in the metric space. The resulting family of lines may have
strange properties. For instance, a line can be a proper subset of
another: in the metric space with points $u,v,x,y,z$ and 
\begin{eqnarray}
\hspace{1.5cm} &&
d(u,v)=d(v,x)=d(x,y)=d(y,z)=d(z,u)=1,\nonumber \\
&& d(u,x)=d(v,y)=d(x,z)=d(y,u)=d(z,v)=2,\nonumber
\end{eqnarray}
we have
\[
\overline{vy}=\{v,x,y\} \;\;\mbox{ and }\;\;
\overline{xy}=\{v,x,y,z\}. 
\]

Chen \cite{Che} proved that a classic theorem of ordered geometry, the
Sylvester-Gallai theorem, generalizes in the framework of metric
spaces (when lines in these spaces are defined differently than
here). Chen and Chv\'{a}tal \cite{CC08} suggested that the de Bruijn -
Erd\H os theorem, too, might generalize in this framework:

\begin{center}
{\em True or false? Every finite metric space $(X,d)$\\
where no line consists of the entire ground set $X$\\
determines at least $|X|$ distinct lines.}
\end{center}
They proved that  
\begin{itemize}

\item in every metric space on $n$ points, there are at least $\lg n$
  distinct lines or else some line consists of all $n$ points.

\end{itemize}
We prove that
\begin{itemize}

\item in every metric space on $n$ points, there are
  $\Omega((n/\rho)^{2/3})$ distinct lines, where $\rho$ is the ratio
  between the largest distance and the smallest nonzero distance
  (Theorem~\ref{range});

\item in every metric space induced by a connected graph on $n$
  vertices, there are $\Omega(n^{2/7})$ distinct lines or else some
  line consists of all $n$ vertices (Corollary~\ref{graphs});

\item in every metric space on $n$ points where each nonzero distance
  equals $1$ or $2$, there are $\Omega(n^{4/3})$ distinct lines and
  this bound is tight (Theorem~\ref{4/3}).

\end{itemize}

\section{Lines in hypergraphs}

A {\em hypergraph\/} is an ordered pair $(X,H)$ such that $X$ is a set
and $H$ is a family of subsets of $X$; elements of $X$ are the {\em
  vertices\/} of the hypergraph and members of $H$ are its {\em
  edges\/}. A hypergraph is called {\em $k$-uniform\/} if each of its
edges consists of $k$ vertices. The definition of lines in a metric
space $(X,d)$ depends only on the $3$-uniform hypergraph $(X,H(d))$
where
\[
H(d)=\{\{a,b,c\}:\; d(a,b)+d(b,c)=d(a,c)
\mbox{ and } a\ne b,\, b\ne c\}:
\]
the line $\overline{uv}$ equals $\{u,v\}\cup\{w:\{u,v,w\}\in H(d)\}.$
This observation suggests extending the notion of lines in metric
spaces to a notion of lines in $3$-uniform hypergraphs: for any two
distinct vertices $u$ and $v$ in a $3$-uniform hypergraph $(X,H)$, the
line $\overline{uv}$ is defined as $\{u,v\}\cup\{w: \{u, v, w\}\in H
\,\}.$ Now every metric space $(X,d)$ and its associated hypergraph
$(X,H(d))$ define the same family of lines.

Let $f(n)$ denote the smallest number of lines in a $3$-uniform
hypergraph on $n$ vertices where no line consists of all $n$ vertices
and let $g(n)$ denote the smallest number of lines in a metric
space on $n$ points where no line consists of all $n$ points. In this
notation, $f(n)\le g(n)$ for all $n$; Chen and Chv\'{a}tal
\cite{CC08} proved that
\[
\lg n \le f(n) < c^{\sqrt{\lg n}}
\]
for some positive constant $c$. (The proof of the lower bound is based
on the observation that $w\not\in\overline{uv}$ if and only if
$v\not\in\overline{uw}$, and so --- unless some line contains all the
vertices --- the mapping that assigns to each vertex the set of lines
containing it is one-to-one.) The upper bound on $f(n)$ does not
rule out the possibility of $g(n)=n$: not all $3$-uniform
hypergraphs arise from metric spaces $(X,d)$ as $(X,H(d))$. (It has
been proved (\cite{Chv,Che}) that the hypergraph consisting of the
seven vertices $0,1,2,3,4,5,6$ and the seven edges $\{ i
\bmod{7},\, (i+1) \bmod{7},\, (i+3) \bmod{7}\}$ with $i=0,1,2,3,4,5,6$
does not arise from any metric space. This $3$-uniform hypergraph is
known as the {\em Fano plane\/} or the {\em projective plane of order
  two\/}.)

We let $K^3_4$ denote the $3$-uniform hypergraph with four vertices
and four edges.

\begin{lem}\label{sparse}
  Let $H$ be a $3$-uniform hypergraph, let $x$ be a vertex of $H$, and
  let $T$ be a set of vertices of $H$ such that {\rm (}i\/{\rm )} $x\not\in T$ and
  {\rm (}ii\/{\rm )} there are no vertices $u,v,w$ in $T$ such that $x,u,v,w$ induce
  a $K^3_4$ in $H$. Then $H$ defines at least $0.25(2|T|)^{2/3}$
  distinct lines.
\end{lem}

\noindent  {\bf Proof.} We may assume that $|T|>4$: otherwise
$0.25(2|T|)^{2/3}\le 1$, which makes the assertion trivial.  Let $S$
denote a largest subset of $T$ such that all the lines $\overline{xv}$ with
$v\in S$ are identical. Now $H$ defines at least $|T|/|S|$ distinct
lines, which gives the desired conclusion when $|S|\le
(2|T|)^{1/3}$. We will prove that $H$ defines at least $|S|(|S|-1)/2$
distinct lines, which gives the desired conclusion when $|S|\ge
(2|T|)^{1/3}\ge 2$. More precisely, we will prove that all the lines
$\overline{uv}$ with $u,v\in S$ are distinct. For this purpose, consider any
three pairwise distinct vertices $u,v,w$ in $S$. Since
$\overline{xu}=\overline{xv}=\overline{xw}$, all three of $\{x,u,v\}$, $\{x,u,w\}$,
$\{x,v,w\}$ are edges of $H$; since $x,u,v,w$ do not induce a $K^3_4$,
it follows that $\{u,v,w\}$ is not an edge of $H$; since $w$ is an
arbitrary vertex in $S$ distinct from $u$ and $v$, the line $\overline{uv}$
intersects $S$ in $\{u,v\}$.  \hfill $\Box$

\section{Lines in metric spaces}

\begin{thm}\label{range}
  In every metric space on $n$ points such that $n\ge 2$, there are at
  least $0.25(n/\rho)^{2/3}$ distinct lines, where $\rho$ is the ratio between
  the largest distance and the smallest nonzero distance.
\end{thm}

\noindent  {\bf Proof.}  Let the metric space be $(X,d)$, let $\delta$
denote the smallest nonzero distance and let $x$ be an arbitrary point
of $X$.  The $n-1$  distances $d(x,u)$ with $u\ne x$ are
distributed into buckets $[i\delta, (i+1)\delta)$ with $i=1,2,\ldots ,
\lfloor \rho\rfloor$. It follows that there are a subset $T$ of
$X-\{x\}$ and a positive integer $i$ such that
\[
u\in T \;\Rightarrow\; i\delta \le d(x,u)< (i+1)\delta
\]
and
\[
|T|\;\ge \; \frac{n-1}{\lfloor \rho\rfloor}\;\ge \; \frac{n-1}{\rho}\;\ge \; \frac{n}{2\rho}\, .
\]
We will complete the proof by showing that the hypergraph $(X,H(d))$
satisfies the hypothesis of Lemma~\ref{sparse}. For this purpose,
consider arbitrary points $u,v,w$ in $T$ such that $\{x,u,v\},
\{x,u,w\}, \{x,v,w\} \in H(d)$; we will prove that $\{u,v,w\}\not\in
H(d)$. Since $\{x,u,v\}\in H(d)$ and $|d(x,u)-d(x,v)|<\delta\le
d(u,v)$, we have $d(u,v)=d(x,u)+d(x,v)$; similarly,
$d(u,w)=d(x,u)+d(x,w)$ and $d(v,w)=d(x,v)+d(x,w)$. Now 
\[
i\delta \le d(x,u), d(x,v), d(x,w) < 2i\delta,
\]
and so 
\[
2i\delta \le d(u,v), d(u,w), d(v,w) < 4i\delta,
\]
and so $\{u,v,w\}\not\in H(d)$.
\hfill $\Box$

\section{Metric spaces induced by graphs}

Every finite connected undirected graph induces a metric space, where
the distance between vertices $u$ and $v$ is defined as the smallest
number of edges in a path from $u$ to $v$.

\begin{thm}\label{diam}
  If, in a metric space $(X,d)$ induced by a graph of diameter $t$, no
  line equals $X$, then there are at least $\sqrt{t/2\;}$ distinct
  lines.
\end{thm}

\noindent  {\bf Proof.}  There are vertices $v_0, v_1, \ldots , v_t$ such
that $d(v_i,v_j)=j-i$ whenever $0\le i<j \le t$. Consider a largest
set $S$ of subscripts $r$ such that all lines $\overline{v_rv_{r+1}}$ are
equal.  There are at least $t/|S|$ distinct lines; this gives the
desired conclusion when $|S|\le \sqrt{2t}$.  We will complete the
argument by proving that there are at least $|S|/2$ distinct lines,
which gives the desired conclusion when $|S|\ge \sqrt{2t}$.

Let $u$ be any vertex outside the line $\overline{v_rv_{r+1}}$ with $r\in
S$. We will prove that at least $|S|/2$ of the lines $\overline{uv_r}$ with
$r\in S$ are pairwise distinct: more precisely, for every three
subscripts $i,j,k$ in $S$, at least two of the three lines
$\overline{uv_i}$, $\overline{uv_j}$, $\overline{uv_k}$ are distinct.

By the triangle inequality and since $u\not\in\overline{v_rv_{r+1}}$, we have 
\[
|\, d(u,v_{r})-d(u,v_{r+1})\, |\;<\; d(v_r,v_{r+1}) \;\;\mbox{ for all $r$ in $S$};
\]
since $d(v_r,v_{r+1})=1$, it follows that $d(u,v_{r})=d(u,v_{r+1})$
for all $r$ in $S$. Now consider any three subscripts $i,j,k$ in $S$
such that $i<j<k$. We have
\begin{eqnarray*}
&d(v_{i},\!v_{j})\!+\!d(v_{j},\!u)>d(v_{i+1},\!v_{j})\!+\!d(v_{j},\!u) \ge d(v_{i+1},\!u) = d(v_{i},\!u),\! &\\
&d(u,\!v_{i})\!+\!d(v_{i},\!v_{j}) > d(u,\!v_{i})\!+\!d(v_{i+1},\!v_{j}) = d(u,\!v_{i+1})\!+\!d(v_{i+1},\!v_{j})\ge d(u,\!v_{j}),\! &
\end{eqnarray*}
and so $v_j\in\overline{uv_i}$ if and only if
$d(v_{i},\!u)\!+\!d(u,\!v_{j})= d(v_{i},\!v_{j})$.  Similarly,
$v_k\in\overline{uv_j}$ if and only if $d(v_{j},\!u)\!+\!d(u,\!v_{k})=
d(v_{j},\!v_{k})$.  Since
\[
d(u,v_i)+d(u,v_{k})\:=\: d(u,v_i)+d(u,v_{k+1})\:\ge \: d(v_i,v_{k+1}) \:=\: k+1-i,
\]
we have $d(u,v_i)>j-i$ or else $d(u,v_{k})>k-j$. If
$d(u,v_i)>j-i$, then $d(v_{i},\!u)\!+\!d(u,\!v_{j}) >
d(v_{i},\!v_{j})$, and so $v_j\not\in\overline{uv_i}$ (and
$v_i\not\in\overline{uv_j}$), which implies $\overline{uv_i}\ne \overline{uv_j}$. If
$d(u,v_{k+1})>k-j$, then $d(v_{j},\!u)\!+\!d(u,\!v_{k}) >
d(v_{j},\!v_{k})$, and so $v_k\not\in\overline{uv_j}$ (and
$v_j\not\in\overline{uv_k}$), which implies $\overline{uv_j}\ne \overline{uv_k}$.
\hfill $\Box$

\begin{cor}\label{graphs}
  If, in a metric space induced by a connected graph on $n$ vertices, no line
  consists of all $n$ vertices, then there are at least
  $2^{-8/7}n^{2/7}$ distinct lines.
\end{cor}

\noindent  {\bf Proof.} If the graph has diameter at most $2^{-9/7}n^{4/7}$,
then the bound follows from Theorem~\ref{range}; else it follows from
Theorem~\ref{diam}. \hfill $\Box$

\bigskip

\section{When each nonzero distance equals $1$ or $2$}

By a {\em $1$-$2$ metric space,\/} we mean a metric space where each nonzero distance is $1$ or $2$. 

\begin{thm}\label{4/3}
  The smallest number $h(n)$ of lines in a $1$-$2$ metric space on $n$ points satisfies 
the inequalities
\[
(1+o(1))\alpha n^{4/3} \le h(n) \le (1+o(1))\beta n^{4/3}
\]
with $\alpha = 2^{-7/3}$ and $\beta =3\cdot 2^{-5/3}$.
\end{thm}

We say that points $u,v$ in a $1$-$2$ metric space are {\em twins\/} if, and only if, $d(u,v)=2$ 
and $d(u,w)=d(v,w)$ for all $w$ distinct from both $u$ and $v$. Our proof of Theorem~\ref{4/3} 
relies on the following lemma, whose proof is routine.

\begin{lem}\label{12triv}
If $u_1,u_2,u_3,u_4$ are four distinct points in a $1$-$2$ metric space, then:
\begin{itemize}
\item [{\rm (i)}] if $d(u_i,u_j)=1$ for all choices of distinct $i$ and $j$, 
then $\overline {u_1u_2}\ne \overline {u_3u_4}$,
\item [{\rm (ii)}] if $d(u_1,u_2)=1$ and $d(u_3,u_4)=2$, 
then $\overline {u_1u_2}\ne \overline {u_3u_4}$,
\item [{\rm (iii)}] if $d(u_1,u_2)=d(u_3,u_4)=2$ and $u_4$ has a twin other than $u_3$,\\
then $\overline {u_1u_2}\ne \overline {u_3u_4}$.
\end{itemize}
If $u_1,u_2,u_3$ are three distinct points in a $1$-$2$ metric space, then:
\begin{itemize}
\item [{\rm (iv)}] if $d(u_1,u_2)=d(u_2,u_3)=1$ and $u_1,u_3$ are not twins, 
then $\overline {u_1u_2}\ne \overline {u_2u_3}$,
\item [{\rm (v)}] if $d(u_1,u_2)=1$, $d(u_2,u_3)=2$, and $u_3$ has a twin other than $u_2$,\\
then $\overline {u_1u_2}\ne \overline {u_2u_3}$,
\item [{\rm (vi)}] if $d(u_1,u_2)=d(u_2,u_3)=2$, 
then $\overline {u_1u_2}\ne \overline {u_2u_3}$.
\end{itemize}
\hfill $\Box$
\end{lem}

\bigskip

\noindent  {\bf Proof of Theorem~\ref{4/3}.} To see that $h(n)\le
(1+o(1))\beta n^{4/3}$, consider the metric space where the ground set
is split into pairwise disjoint groups of sizes as nearly equal as
possible, every two points that belong to two different groups have
distance $1$, and every two points that belong to one group have
distance $2$. If each group includes at least three points, then
$\overline{uv}=\overline{wx}$ if and only if either $\{u,v\}=\{w,x\}$ or else
there are two distinct groups such that each of the sets $\{u,v\}$,
$\{w,x\}$ has one element in each of these two groups.  Consequently,
when there are $n$ points altogether and $(1+o(1))2^{-1/3}n^{2/3}$
groups, there are $(1+o(1))\beta n^{4/3}$ lines.

To prove that $h(n)\ge (1+o(1))\alpha n^{4/3}$, consider an arbitrary
$1$-$2$ metric space $(X,d)$ and write $n=|X|$. Let $X_1$ be any
maximal subset of $X$ that does not contain a pair of twins.

{\sc Case 1:} $|X_1|\ge n/2$. 
In this case, consider a largest set of distinct two-point subsets $\{u_i,v_i\}$ 
($i=1,2,\ldots ,s$) of $X_1$ such that
\[
\overline{u_{1}v_{1}} = \overline{u_{2}v_{2}} = \ldots = \overline{u_{s}v_{s}}.
\]
Since every two-point subset of $X_1$ determines a line, there are at least
\[
\binom{|X_1|}{2}\cdot\frac{1}{s}
\]
distinct lines; this gives the desired conclusion when $s\le
(n/2)^{2/3}$. We will complete the argument by proving that there are
at least
\[
\binom{s}{2} - 5
\]
distinct lines, which gives the desired conclusion when $s\ge (n/2)^{2/3}$. 

For this purpose, we may assume that $s\ge 5$. Part (iv) of
Lemma~\ref{12triv} guarantees that the sets $\{u_i,v_i\}$ with
$d(u_i,v_i)=1$ are pairwise disjoint; part (iv) of Lemma~\ref{12triv}
guarantees that the sets $\{u_i,v_i\}$ with $d(u_i,v_i)=2$ are
pairwise disjoint; part (ii) of Lemma~\ref{12triv} guarantees that
each of the sets $\{u_i,v_i\}$ with $d(u_i,v_i)=1$ meets each of the
sets $\{u_i,v_i\}$ with $d(u_i,v_i)=2$; now our assumption $s\ge 5$ guarantees that 
all $s$ distances $d(u_i,v_i)$ are equal.
We are going to prove that there are at least $s(s-1)/2$ distinct
lines: for every choice of subscripts $i,j$ such that $1\le i < j \le
s$, there is a line $L_{ij}$ such that
\[
\{u_k,v_k\}\subseteq L_{ij} \;\;\Leftrightarrow\;\; k\in \{i,j\}.
\]

{\sc Subcase 1.1:} $d(u_1,v_1) = d(u_2,v_2) = \ldots = d(u_s,v_s) =1$.

Since $\{u_j,v_j\}\subseteq \overline{u_{j}v_{j}}=\overline{u_{i}v_{i}}$ and
$\{u_i,v_i\}\subseteq \overline{u_{i}v_{i}}=\overline{u_{j}v_{j}}$, we may assume
(after switching $u_j$ with $v_j$ if necessary) that $d(u_i,u_j) = 2$
and $d(u_i,v_j) = d(u_j,v_i) = 1$. Now we may set
$L_{ij}=\overline{u_{i}u_{j}}$: if $k\not\in \{i,j\}$, then $u_j\in
\overline{u_{j}v_{j}}=\overline{u_{k}v_{k}}$ implies that one of $d(u_j,u_k)$ and
$d(u_j,v_k)$ equals $2$, and so $\{u_k,v_k\}\not\subseteq
\overline{u_{i}u_{j}}$.

{\sc Subcase 1.2:} $d(u_1,v_1) = d(u_2,v_2) = \ldots = d(u_s,v_s) =2$.

Since $\overline{u_{1}v_{1}} = \overline{u_{2}v_{2}} = \ldots = \overline{u_{s}v_{s}}$,
the distance between any point in one of the sets $\{u_1,v_1\}$,
$\{u_2,v_2\}$, \ldots , $\{u_s,v_s\}$ and any point in another of
these $s$ sets equals $1$; it follows that we may set
$L_{ij}=\overline{u_{i}u_{j}}$.

{\sc Case 2:} $|X_1|< n/2$. 
Write $X_2=X-X_1$, consider a largest set $S$ 
of points in $X_2$ such that 
\[
u,v\in S,\, u\ne v\;\;\Rightarrow\;\; d(u,v)=1,
\]
and write
\begin{eqnarray*}
E_1 &=& \{\;\{u,v\}\!: u,v\in S,\, u\ne v \},\\
E_2 &=& \{\;\{u,v\}\!: u,v\in X_2,\, d(u,v)=2 \}.
\end{eqnarray*}
Since every vertex of $X_2$ has a twin (else it could be added to
$X_1$), Lemma~\ref{12triv} guarantees that every two distinct pairs in
$E_1\cup E_2$ determine two distinct lines.  We complete the argument
by pointing out that
\[
|E_1\cup E_2|\ge (1+o(1))\alpha n^{4/3}:
\]
the famous theorem of Tur\' an~\cite{T} guarantees that
\[
|S| \;\ge \; \frac{|X_2|^2}{2|E_2|+|X_2|}, 
\]
and so $|E_2|\;<\; \alpha n^{4/3}$ implies $|E_1| \;\ge \;
(1+o(1))\,\alpha n^{4/3}$. \hfill $\Box$

\bigskip

The lower bound of Theorem~\ref{4/3} can be easily improved through a
more careful analysis of Case 2: a routine exercise in calculus shows 
that
\[
x\ge 3,\, y\ge 0 \;\;\Rightarrow\;\;
\frac{1}{2}\cdot \left(\frac{x^2}{2y+x}\right)^2 + y \;\ge \; \beta x^{4/3}-\frac{x}{2}\: , 
\]
and so $|E_1\cup E_2|\ge (1+o(1))\beta |X_2|^{4/3}$. Perhaps 
$h(n)=(1+o(1))\beta n^{4/3}$.   

\bigskip

\begin{center}
{\bf Acknowledgment}
\end{center}

\noindent 
This research was carried out in ConCoCO (Concordia Computational
Combinatorial Optimization Laboratory) and undertaken, in part, thanks
to funding from the Canada Research Chairs Program and from the
Natural Sciences and Engineering Research Council of Canada.

\end{document}